\documentclass{amsart}

\numberwithin{equation}{section}

\theoremstyle{plain}
\newtheorem{theorem}{Theorem}[section]
\newtheorem{THEOREM}{Theorem}

\newtheorem{proposition}[theorem]{Proposition}
\newtheorem{lemma}[theorem]{Lemma}
\newtheorem{corollary}[theorem]{Corollary}

\theoremstyle{definition}
\newtheorem{definition}[theorem]{Definition}
\newtheorem*{definition*}{Definition}

\newtheorem*{notation*}{Notation}

\theoremstyle{remark}
\newtheorem{remark}[theorem]{Remark}
\newtheorem*{remark*}{Remark}

\newtheorem*{remarks*}{Remarks}
\newtheorem{example}[theorem]{Example}
\newtheorem*{example*}{Example}

\newcounter{varenum}

\newcommand\st{\vert}						
\newcommand\eqd{:=}						
\newcommand\Z{\mathbb{Z}}					
\newcommand\C{\mathbb{C}}					
\newcommand\GL{{\mathrm{GL}}}					
\newcommand\SL{{\mathrm{SL}}}					
\newcommand\gl{{\mathfrak{gl}}}				
\newcommand\gal{\mathop{\mathrm{Gal}}\nolimits}		
\newcommand\Pic{\mathop{\mathbf{Pic}}\nolimits}		
\newcommand\limind{\varinjlim}				
\newcommand\spec{\mathop{\mathrm{Spec}}\nolimits}	
\newcommand\spf{\mathop{\mathrm{Spf}}\nolimits}		
\newcommand\res{\mathop{\mathrm{res}}}			
\newcommand\id{{\mathrm{id}}}					
\newcommand\gm{\mathbf{G_m}}					
\newcommand\iso{{\widetilde\to}}				
\newcommand\tr{\mathop{\mathrm{tr}}\nolimits}		
\newcommand\ad{\mathop{\mathrm{ad}}\nolimits}		
\newcommand\coker{\mathop{\mathrm{coker}}\nolimits}	
\newcommand\p[1]{{\mathbb{P}^{\mathbf{#1}}}}		

\newcommand\tA{{\widetilde A}}
\newcommand\tX{{\widetilde X}}
\newcommand\tx{{\tilde x}}
\newcommand\tz{{\tilde z}}
\newcommand\cG{{\mathcal{G}}}
\newcommand\cF{{\mathcal{F}}}

\newcommand\gfr{{\mathfrak{g}}}
\newcommand\hfr{{\mathfrak{h}}}
\newcommand\tL{{\widetilde L}}
\newcommand\tnabla{{\widetilde\nabla}}

\newcommand\SCurv{{\mathrm{SCurv}}}
\newcommand\MS{{{\mathbf M}_\sharp}}
\newcommand\cMS{{M_\sharp}}
\newcommand\pS{{p_\sharp}}
\newcommand\pH{{p_H}}
\newcommand\al{{\bullet}}
\newcommand\tXS{{\tX_\sharp}}
\newcommand\lS{{l_\sharp}}
\newcommand\partialS{{\partial_\sharp}}
\newcommand\Connst{\mathbf{Conn}}
\newcommand\Higgst{\mathbf{Higgs}}
\newcommand\Cn{{{\Connst}_\lambda}}		
\newcommand\Cnu{{\Connst}}			
\newcommand\Hg{{\Higgst}} 			
\newcommand\Hgp{{\Hg'}}				
\newcommand\Cnf{{{\Connst}_{form}}}		
\newcommand\Cnfp{{{\Connst}'_{form}}}	
\newcommand\Conngr{\mathcal{C}onn}
\newcommand\Higgsgr{\mathcal{H}iggs} 
\newcommand\Conn[1]{{\Conngr_\lambda(#1)}}
\newcommand\tConn[1]{{{\widetilde{\Conngr}}_\lambda(#1)}}
\newcommand\Higgs[1]{{\Higgsgr(#1)}}
\newcommand\tHiggs[1]{{\widetilde{\Higgsgr}(#1)}}
\newcommand\CC[1]{{\mathcal{C}(#1)}}
\newcommand\tCC[1]{{\widetilde{\mathcal{C}}(#1)}}
\newcommand\F{\mathcal F}
\newcommand\Vect{{\mathcal{V}ect}}		
\newcommand\Princ{{\mathcal{P}rinc}}	

\newcommand\cnt{\mathfrak{c}}
\newcommand\ccnt{\mathfrak{c}^\vee}

\title{Moduli of connections with a small parameter on a curve}
\author{D.~Arinkin}
\address{Department of Mathematics, University of Chicago, Chicago, Illinois 60637, USA}
\email{arinkin@math.uchicago.edu}
\thanks{Partially supported by NSF grants DMS-0100108 and DMS-0401164}

\begin{document}

\begin{abstract}
We study $\GL_2$-bundles with connections with a small parameter on a smooth projective curve. We describe
an open subset in the moduli space of such bundles. The description degenerates into
the Hitchin fibration as the parameter tends to zero.
\end{abstract}

\maketitle

\section{Introduction}

\subsection{}
The moduli space of Higgs bundles on a curve admits a well-known description in terms of spectral curves 
(the Hitchin fibration). On the other hand, Higgs bundles can be viewed as a degeneration of bundles with 
connections: P.~Deligne introduced the notion of `$\lambda$-connections', and Higgs fields (resp. connections) 
are $\lambda$-connections for $\lambda=0$ (resp. $\lambda=1$).
It is natural to ask whether spectral curves can be used to describe the moduli space of $\lambda$-connections 
for $\lambda\ne0$.

The simplest case is when $\lambda\in\C[[\lambda]]$ is a formal parameter; that is, the $\lambda$-connections 
considered are formal deformations of Higgs bundles. This case has the following advantage:
if a $\lambda$-connection is a formal deformation of a Higgs bundle, we can try using the spectral curve corresponding
to the Higgs bundle to describe the $\lambda$-connection. Informally, if $\lambda$ is an actual number
rather then a formal parameter (for instance, $\lambda=1$), we would not know which
spectral curve to use.  

Let $\Cn$ be the moduli space of $\lambda$-connections; more precisely, it parametrizes triples $(L,\nabla,\lambda)$,
where $L$ is a $G$-bundle on $X$, $\lambda\in\C$, and $\nabla$ is a $\lambda$-connection on $L$. Here $X$ is a smooth
curve and $G$ is a reductive group. The moduli stack of Higgs bundles is the closed substack
$\Hg\subset\Cn$ given by the condition $\lambda=0$. Making $\lambda$ a formal parameter corresponds to working with
the formal completion $\Cnf$ of $\Cn$ along $\Hg$ instead of $\Cn$ itself.

\begin{remark}
Although $\lambda$-connections are interesting geometric objects in their own right, they are particularly important
because they can be used to compactify the moduli stack $\Cnu$ of ordinary connections (\cite{Si1}, \cite{Si2}). 
We hope that studying $\Cnf$ can improve our understanding of $\Cnu$ (which is important, for instance, in the
geometric Langlands program). One particular case ($\SL_2$-bundles
with connections on $\p1$ with four simple poles) appears in \cite{Ar}: a statement about $\Cnf$ 
(\cite[Proposition 6]{Ar}) is used to compute the cohomology groups 
$H^i(\Cnu,F)$ for some natural coherent sheaves $F$. 
\end{remark}

The problem simplifies further if we consider formal deformations of only those Higgs bundles that are non-degenerate
in some sense. Geometrically, this corresponds to taking an open substack $\Hgp\subset\Hg$ (parametrizing non-degenerate
bundles) and studying  the formal completion $\Cnfp$ of $\Cn$ along $\Hg$. Let us look at three different non-degeneracy
conditions for Higgs bundles.

Firstly, let us consider only the Higgs bundles with unramified spectral curves (equivalently,
the Higgs field is regular semisimple). The answer in this case is relatively simple (see Theorem \ref{th:Wasow})
and goes back to W.~Wasow. However, from the geometric point of view, this non-degeneracy condition is too restrictive.
Geometrically, the most interesting situation is when $X$ is projective (so that $\Hg$ and $\Cn$
are algebraic stacks); but if $X$ is projective and not elliptic,
it has no unramified spectral curves.

Secondly, let us allow ramifications, but only of `the simplest possible kind'. More precisely, if $G=\GL_2$, we 
consider Higgs bundles whose spectral curves are smooth
(possibly ramified) covers of $X$. For arbitrary reductive $G$, roughly speaking, we consider 
Higgs bundles that can be locally on $X$ reduced to Higgs bundles over $\GL_2$ with smooth spectral curves (see Remark
\ref{rm:cameral} for the precise condition). The main result of this paper describes $\Cnfp$ via spectral curves for this 
non-degeneracy condition.

\begin{remark} For arbitrary reductive $G$, this non-degeneracy condition is most easily formulated using the notion
of cameral covers. Recall (\cite{DG}) that a cameral cover is a cover of $X$ that is locally isomorphic to the pull-back
of the \emph{universal cameral cover} $\hfr\to\hfr/W$, where $\hfr$ is the Cartan algebra of $G$ and $W$ is the Weyl group.
To a Higgs bundle on $X$, there corresponds a cameral cover $X_{cam}\to X$; locally on $X$, a Higgs field is essentially a map 
$X\to\gfr$ (where $\gfr$ is the Lie algebra of $G$), and $X_{cam}\to X$ is the pull-back of the universal cameral cover 
under the composition
\begin{equation*}
X\to\gfr\to\gfr/G=\hfr/W.
\end{equation*}
The second non-degeneracy condition on a Higgs bundle is that $X_{cam}$ is smooth. 
\label{rm:cameral}
\end{remark}

Finally, the last non-degeneracy condition allows a more general kind of ramifications. For instance, if $G=\GL_m$,
the spectral curve is a degree $m$ cover of $X$. Then let us work with Higgs bundles whose spectral curves are smooth.
If $m>2$, then this is a more general (and more natural) condition than the previous one, which only allows smooth
curves whose ramification points have degree 2. It is possible to define this last non-degeneracy condition for an arbitrary
reductive group, not just for $\GL_m$. However, to use this condition one needs to work with non-smooth cameral covers,
which is more complicated. Moduli of $\lambda$-connections for this non-degeneracy conditions will be studied elsewhere.

In this paper, we work with the second non-degeneracy condition. We set $G=\GL_2$ and
consider Higgs bundles whose spectral curves are smooth (the results can then be extended to other reductive groups
by using Levi subgroups). We use spectral curves to describe $\lambda$-connections 
that are formal deformations of such Higgs bundles (Theorem \ref{th:lambda-connections}), and then derive a
description of $\Cnfp$ if $X$ is projective (Theorem \ref{th:lambda-connections via foliations}). 

\subsection{Conventions and notation}
In this work, the ground field is $\C$, that is, `scheme' means `$\C$-scheme', $\GL_2$ means $\GL_2(\C)$, and
so on. However, our methods are purely algebraic, so our results hold over any algebraically closed field
of characteristic zero.

For a scheme (or a formal scheme, or a stack) $S$ and an integer $m>0$, consider the following three versions
of the category of $\GL_m$-bundles on $S$:
\begin{itemize}
\item the category $\Vect$ of rank $m$ vector bundles on $S$;
\item the subcategory $\Vect^\times\subset\Vect$ of `invertible arrows': the objects of $\Vect^\times$ are
rank $m$ vector bundle of $S$, but arrows $L_1\to L_2$ are isomorphisms between
$L_1$ and $L_2$, rather than all homomorphisms ($L_1,L_2\in\Vect^\times$);
\item the category $\Princ$ of principal $\GL_m$-bundles on $S$. 
\end{itemize}
The categories $Vect^\times$ and $\Princ$ are naturally equivalent. $\Vect$ is a pre-additive category 
(morphisms between any two objects form a vector space), while $\Princ$ and $\Vect^\times$ are groupoids 
(all morphisms are invertible). Note also that $\Princ$ makes sense for groups other than $\GL_m$, but
$\Vect$ does not.

In this paper, we work with $\Vect^\times$, which we call \emph{the groupoid of $\GL_m$-bundles on $S$} (or
\emph{$\gm$-bundles}, if $m=1$). The same convention applies to $\GL_m$-bundles with additional structures, such as connections,
$\lambda$-connections, or Higgs fields. It is interesting to note however that Theorem \ref{th:lambda-connections}
also holds in the `pre-additive' settings; the proof is left to the reader.

\section{Main results}

\subsection{$\lambda$-connections on a curve}

Although our results hold for arbitrary reductive group, we prefer to formulate them for $\GL_2$. Let $X$ be a smooth 
projective curve over $\C$.

\begin{definition} Let $L$ be a $\GL_m$-bundle on $X$; that is, $L$ is a locally free sheaf of rank $m$. 
A \emph{$\lambda$-connection} on $L$ (for some $\lambda\in\C$) is a $\C$-linear map $\nabla:L\to L\otimes\Omega_X$
which satisfies the $\lambda$-Leibniz identity:
\begin{equation}
\nabla(fs)=f\nabla(s)+\lambda s\otimes df
\label{lambda-Leibniz}
\end{equation}
for any $f\in O_X$, $s\in L$. 
\end{definition}

\begin{example} If $\lambda=1$, we get the usual Leibniz identity, so $\nabla$ is a connection
on $L$. More generally, for any $\lambda\ne0$, $\lambda^{-1}\nabla$ is a connection
on $L$. On the other hand, if $\lambda=0$, $\nabla$ is $O_X$-linear (a \emph{Higgs field}).
\end{example}

Denote by $\Hg$ the moduli stack of \emph{Higgs bundles} $(L,\nabla)$ over $\GL_2$ on $X$; that
is, $L$ is a $\GL_2$ bundle on $X$ and $\nabla$ is a $\lambda$-connection on $L$ for $\lambda=0$.
The stack $\Hg$ has a well-known geometric description (the Hitchin
fibration) via spectral curves, which we remind in Theorem \ref{th:Hitchin}. Our aim is to provide a similar description 
for the moduli stack of $\lambda$-connections 
(or at least its open subset) when $\lambda\in\C[[\lambda]]$ is a formal parameter.

Denote by $\Cn$ the moduli stack of triples $(L,\nabla,\lambda)$, where $\lambda\in\C$,
$L$ is a $\GL_2$-bundle on $X$ and $\nabla$ is a $\lambda$-connection on $L$. Then $\Hg$ is identified
with the closed substack of $\Cn$ formed by triples $(L,\nabla,\lambda)$ with $\lambda=0$. 
We will describe an open subset in the formal completion of $\Cn$ along $\Hg$.

Let us recall the geometric description of $\Hg$ (\cite{Hi}, see also \cite{DG} for a much 
more general statement). 
Let $p:T^*X\to X$ be the cotangent bundle. 

\begin{definition} A pure dimension $1$ subscheme $\tX\subset T^*X$ is a \emph{spectral curve} (for $\GL_2$)
if the projection $p_\tX\eqd p|_\tX:\tX\to X$ is finite of degree 2. 
\end{definition}

Let $\mu=\mu_\tX\in H^0(\tX,p_\tX^*\Omega_X)$ 
be the restriction of the natural 1-form $\mu_{T^*X}\in H^0(T^*X,p^*\Omega_{X})$. Denote by $\SCurv$ the space of all
spectral curves. $\SCurv$ is isomorphic to an affine space: the coordinates on 
$\SCurv$ are the coefficients of the equation for $\tX$.

\begin{theorem} Let $(L,\nabla)\in\Hg$ be a Higgs bundle.

\begin{enumerate}
\item There exists a unique spectral curve $\tX\in\SCurv$ and a unique (up to a canonical isomorphism) 
coherent $O_{\tX}$-module $l$ such that
$L=(p_\tX)_*l$ and  $\nabla=(p_\tX)_*\mu$. We call $(\tX,l)$ the \emph{spectral data} of $(L,\nabla)$.

\item If $\tX$ is smooth, $l$ is an invertible sheaf on $\tX$.

\item For a smooth spectral curve $\tX$ and an invertible sheaf $l$ on $\tX$, there is a unique
(up to a canonical isomorphism) $(L,\nabla)\in\Hg$ such that
$(\tX,l)$ is the spectral data of $(L,\nabla)$. 
\end{enumerate}
\qed
\label{th:Hitchin}
\end{theorem}

\begin{remark} Consider the morphism $\pH:\Hg\to\SCurv$ that sends a Higgs bundle to its spectral curve
(the Hitchin fibration).
Then Theorem \ref{th:Hitchin} implies that the fiber of $\pH$ over a smooth spectral curve $\tX\in\SCurv$ is
the moduli stack $\Pic(\tX)$ of line bundles on $\tX$. 
\label{rm:Hitchin}
\end{remark}

Our first result is a version of Theorem \ref{th:Hitchin} for $\lambda$-connections.
Let us start with some definitions.

\begin{definition} \emph{A $\C[[\lambda]]$-family} of $\GL_m$-bundles with $\lambda$-connections on a smooth curve 
$X$ is a pair $(L,\nabla)$, where $L$ is a $\GL_m$-bundle on the formal
scheme $X[[\lambda]]\eqd \limind X\times\spec\C[\lambda]/(\lambda^i)$, and 
$\nabla:L\to L\otimes_{O_X}\Omega_X$ is a $\C[[\lambda]]$-linear $\lambda$-connection. 

The \emph{reduction} of $(L,\nabla)$ modulo $\lambda$ is the Higgs bundle $(L_0,\nabla_0)$ on $X$, 
where $L_0\eqd L/\lambda L$ is the restriction of $L$ to $X$, and $\nabla_0:L_0\to L_0\otimes\Omega_X$ is induced
by $\nabla$. 
\label{df:C[[lambda]]-family}
\end{definition}

Let us fix a smooth (but not necessarily projective) curve $X$ and a smooth spectral curve $\tX\subset T^*X$.
Denote by $\Conn{\tX}$ the groupoid of $\C[[\lambda]]$-families $(L,\nabla)$ of $\GL_2$-bundles with connections 
on $X$ such that $\tX$ equals the spectral curve of $(L_0,\nabla_0)$ (where $(L_0,\nabla_0)$ is the reduction
of $(L,\nabla)$ modulo $\lambda$). We would like to describe $\Conn{\tX}$ in terms of bundles on the spectral curve $\tX$.

Denote by $\tConn{\tX}$ the groupoid of $\C[[\lambda]]$-families $(l,\delta)$ of $\gm$-bundles with
connections on $\tX$ such that $\delta:l\to l\otimes\Omega_\tX(\tx_1+\dots+\tx_n)$ has first
order poles at $\tx_1,\dots,\tx_n$ (the ramification locus of $p_\tX:\tX\to X$), the residue of $\delta$ at $\tx_i$
equals $-\lambda/2$ (the notion of residue of a $\lambda$-connection is straightforward), and the reduction 
$\delta_0$ of $\delta$ modulo $\lambda$ equals $\mu\in H^0(\tX,\Omega_\tX)$. Notice that $\delta_0$ is a Higgs field
on the line bundle $l/\lambda l$, and a Higgs field on a line bundle is just a differential form.

\begin{THEOREM} \label{th:lambda-connections}
There exists an equivalence of categories $\F:\Conn{\tX}\iso\tConn{\tX}$ for a smooth 
(not necessarily projective) curve $X$ and a smooth spectral curve $\tX\subset T^*X$.
\end{THEOREM}

\begin{remark*} 
Theorem \ref{th:lambda-connections} is significantly simplified if $\tX$ is unramified over $X$ 
(see Theorem \ref{th:Wasow}). This special case goes back to Wasow (\cite[Theorem 25.2]{Wa}). 
\end{remark*}
\begin{remark}
\label{rm:lambda-connections on line bundles}
The groupoid $\tConn{\tX}$ has a simpler description. Namely, for any $(l,\delta)\in\tConn{\tX}$, the formula
$\partial\eqd \lambda^{-1}(\delta-\mu)$ defines a connection $\partial:l\to l\otimes\Omega_\tX(\tx_1+\dots+\tx_n)$. In this
manner, $\tConn{\tX}$ identifies with the groupoid of pairs $(l,\partial)$, where $l$ is a line bundle on $\tX[[\lambda]]$,
and $\partial:l\to l\otimes\Omega_\tX(\tx_1+\dots+\tx_n)$ is a ($\C[[\lambda]]$-linear) connection whose residues
at $\tx_i\in\tX$ equal $-1/2$.
\end{remark}
The formulation of Theorem \ref{th:lambda-connections} is somewhat unsatisfactory, because the equivalence $\F$ is not described.
However, there are some natural compatibility conditions on $\F$. For instance, if
$X'\subset X$ is an open set, $\tX'\eqd X'\times_X\tX$ is a spectral curve over $X'$, and it is natural to ask
that $\F$ commutes with the restriction functors $\Conn{\tX}\to\Conn{\tX'}$, $\tConn{\tX}\to\tConn{\tX'}$;
essentially, this corresponds to viewing $\Conn{\tX}$, $\tConn{\tX}$ as stacks in the Zariski topology (the \'etale
topology also works). Also, one naturally wants $\F$ to be compatible with Theorem \ref{th:Hitchin}: for
$(L,\nabla)\in\Conn{\tX}$, the spectral data of its reduction $(L_0,\nabla_0)$ should be 
canonically isomorphic to $(\tX,l/\lambda l)$,
where $(l,\delta)=\F(L,\nabla)$. In some sense, the compatibility conditions determine $\F$ up to a unique isomorphism,
see Theorem \ref{th:lambda-connections2}.

\subsection{Moduli space of $\lambda$-connections}

We saw that Theorem \ref{th:Hitchin} provides a geometric description of an open substack of $\Hg$
(Remark \ref{rm:Hitchin}). Similarly, Theorem \ref{th:lambda-connections} can be used to
describe an open stack in the completion of $\Cn$ along $\Hg$. 

Let $X$ be a smooth projective curve. Denote by $\MS$ the moduli stack of collections $(\tX,l,\partial)$, where $\tX\in\SCurv$ is a smooth
spectral curve, $l$ is a line bundle on $\tX$, and $\partial:l\to l\otimes\Omega_\tX(\tx_1+\dots+\tx_n)$ 
is a connection (not a $\lambda$-connection) whose residues at $\tx_1,\dots,\tx_n$ equal $-1/2$. As before,
$\tx_1,\dots,\tx_n$ are the ramification points of $p_\tX:\tX\to X$. 

Consider the projection 
\begin{equation*}
\pS:\MS\to\Hg:(\tX,l,\delta)\mapsto(\tX,l);
\end{equation*} 
here we use Theorem \ref{th:Hitchin} to identify Higgs bundles with their spectral data $(\tX,l)$. The fiber
of $\pS$ over $(\tX,l)$ is the space of connections $\partial:l\to l\otimes\Omega_\tX(\tx_1+\dots+\tx_n)$,
$\res_{\tx_i}\partial=-1/2$. The following statement is immediate:

\begin{lemma} Denote by $\Hgp\subset\Hg$ the open substack of Higgs 
bundles whose spectral data $(\tX,l)$ satisfy two conditions: $\tX$ is smooth, and $\deg(l)=n/2=2g-2$, 
where $n$ is the number of ramification points of $p_\tX$ and $g$ is the genus of $X$. Equivalently,
$(L,\nabla)\in\Hgp$ if its spectral curve is smooth and $\deg(L)=0$. Then

\begin{enumerate}
\item  $\pS(\MS)=\Hgp$.

\item  For $(\tX,l)\in\Hgp$, the fiber $p_2^{-1}(\tX,l)$ is an affine space; the corresponding vector space
is $H^0(\tX,\Omega_\tX)$. More precisely: as $\tX$ varies, the spaces $H^0(\tX,\Omega_\tX)$ form a vector bundle on $\Hgp$,
and $p_2:\MS\to\Hgp$ is a torsor over this vector bundle.
\end{enumerate}
\qed
\end{lemma}

Denote by $\zeta_0$ the relative tangent bundle to $\pS$; it is a foliation on $\MS$, and $\Hgp$ can be viewed
as the quotient of $\MS$ modulo $\zeta_0$. 

\begin{remark} Technically, $\MS$ is an algebraic stack rather then a scheme, and the notion of a foliation
on a stack requires clarification. However, the stack structure on $\MS$ (and on $\Hgp$) is rather
simple: the automorphism group of every point equals $\gm$; that is, $\MS$ is a $\gm$-gerbe over the corresponding
coarse moduli space, $\cMS$. If we choose to work with $\cMS$ instead of $\MS$, then $\zeta_0$ becomes
just a foliation on a smooth algebraic space; the downside is that in this way 
we get a description of the coarse moduli space of $\lambda$-connections rather then the true moduli stack.
We could avoid this difficulty if we rigidify the moduli problem, for instance, by adding a framing of vector
bundles at some points. 

On the other hand, it is not hard to define the notion of a foliation on an algebraic stack (for instance,
using Lee algebroids).
From now on, we will ignore this difficulty and freely use foliations on $\MS$.
\end{remark} 

Notice that $\MS$ carries another foliation, which is defined via isomonodromic deformation. Let us consider the
composition $\pH\circ\pS:\MS\to\SCurv$. The fiber of $\pH\circ\pS$ over a smooth spectral curve $\tX\in\SCurv$ is 
canonically
identified with fibers over infinitesimally close spectral curves (the fiber is essentially the space of rank 1 local
systems on $\tX$ with monodromy $-1$ around the ramification points; therefore, the fiber does not change under 
deformations of $\tX$). More precisely, the morphism $\pH\circ\pS:\MS\to\SCurv$ carries a connection. Let 
$\zeta_\infty$ be the foliation (on $\MS$) of horizontal vector fields with respect to this connection.

Let us now consider $\zeta_0$ and $\zeta_\infty$ as abstract vector bundles (rather then foliations) on $\MS$. 
Over a point $(\tX,l,\partial)\in\MS$, the fiber
of $\zeta_0$ equals $H^0(\tX,\Omega_\tX)$, while the fiber of $\zeta_\infty$ equals $H^0(\tX,N_\tX)$, where $N_\tX$
is the normal bundle to $\tX\subset T^*X$. The symplectic structure on $T^*X$ identifies $N_\tX$ with $\Omega_\tX$;
therefore, $\zeta_0$ and $\zeta_\infty$ are isomorphic as vector bundles on $\MS$.

\begin{remark} For the isomorphism $\Omega_\tX\iso N_\tX$, there are two choices that differ by sign; we choose
the sign so that the diagram
\begin{equation*}
\begin{array}{ccc}
p_{\tX}^*\Omega_X&\to&T(T^*X)|_\tX\cr
\downarrow&&\downarrow\cr
\Omega_\tX&\iso&N_\tX
\end{array}
\end{equation*}
commutes. Here $T(T^*X)|_\tX$ is the restriction to $\tX\subset T^*X$ of the tangent bundle to $T^*X$, the map
$p_{\tX}^*\Omega_X\to T(T^*X)|_\tX$ identifies $p_{\tX}^*\Omega$ with the subbundle of vertical vector fields,
$p_{\tX}^*\Omega_X\to\Omega_\tX$ is the pull-back map for differential forms, and $T(T^*X)|_\tX\to N_\tX$
is the natural projection.
\label{rm:isobundles}
\end{remark}

\begin{definition} Let $\zeta_0,\zeta_\infty\subset TM$ be distributions on a smooth variety $M$, and let
$\nu:\zeta_0\iso\zeta_\infty$ be an isomorphism of vector bundles on $M$. The \emph{linear combination}
$\alpha\zeta_0+\beta\zeta_\infty\subset TM$ is the distribution on $M$ that is the
image of the morphism $\alpha(\id_{\zeta_0})+\beta\nu:\zeta_0\to TM$,
provided the morphism is an embedding of vector bundles. Clearly, the linear combination 
$\alpha\zeta_0+\beta\zeta_\infty$, if it exists, depends only on the ratio $(\alpha:\beta)\in\p1$.

Notice that a linear combination is not necessarily a foliation even if $\zeta_0$ and $\zeta_\infty$ are foliations. 
\end{definition}

\begin{THEOREM} Let $\MS$, $\zeta_0$, and $\zeta_\infty$ be as above, and let us use the isomorphism 
$\zeta_0\iso\zeta_\infty$ from Remark \ref{rm:isobundles} to construct the linear combination
$\zeta_\lambda\eqd\zeta_0-\lambda\zeta_\infty$, $\lambda\in\C$.
\begin{enumerate} 
\item $\zeta_\lambda=\zeta_0-\lambda\zeta_\infty$ is a foliation on $\MS$ for any $\lambda\in\C$.
\item The quotient $\MS/\zeta_\lambda$ exists if $\lambda\in\C[[\lambda]]$ is a formal parameter, and such quotients
form a family $\MS[[\lambda]]/\zeta_\lambda\to\spf\C[[\lambda]]$ over the formal disc.
\item $\MS[[\lambda]]/\zeta_\lambda$ is canonically isomorphic to the formal completion of $\Cn$ along $\Hgp$. This 
isomorphism respects the projection to $\spf\C[[\lambda]]$ (intuitively, $\MS/\zeta_\lambda$ is identified
with an open substack in the moduli stack of $\lambda$-connections when $\lambda$ is a formal parameter).
\end{enumerate}
\label{th:lambda-connections via foliations}
\end{THEOREM} 

\begin{remark} 
\label{re:combination exists}
Let us show that the linear combination $\zeta_\lambda$ exists (as a distribution) for any $\lambda\in\C$. Indeed,
if $\lambda=0$, then $\zeta_\lambda=\zeta_0$; so it is enough to analyze the case $\lambda\ne 0$. Now consider the 
projection $\pH\circ\pS:\MS\to\SCurv$. Its differential $d(\pH\circ\pS)$ vanishes on
$\zeta_0$ and induces an isomorphism between $\zeta_\infty\subset T\MS$ and the pull-back of the tangent bundle from
$\SCurv$ to $\MS$. Therefore, $d(\pH\circ\pS)$ also induces an
isomorphism between $\zeta_\lambda$ and this pull-back. This implies the statement.
\end{remark}

\begin{remark} Note that although $\Hgp$ is open in $\Hg$, it is not dense. Actually, $\Hg$ is disconnected;
its connected components are 
\begin{equation*}
\Hg^{(k)}\eqd\{(L,\nabla)\in\Hg\st\deg(L)=k\}\quad (k\in\Z),
\end{equation*} 
and $\Hgp\subset\Hg^{(0)}$.
However, only the neighborhood of $\Hg^{(0)}\subset\Cn$ is interesting, because $\Hg^{(k)}\subset\Cn$ is a 
connected component
of $\Cn$ for $k\ne0$. This is easy to see by using the exterior square (the `trace') of $\lambda$-connections.
\end{remark}

\subsection{Organization} In Section \ref{sc:reduction}, we formulate a more precise version of 
Theorem \ref{th:lambda-connections}. We then show that the theorem follows from its `formal' version (Theorem 
\ref{th:lambda-connectionsf}), in which $X$ is a formal disc rather than a curve. Theorem \ref{th:lambda-connectionsf}
is proved in Section \ref{sc:proof}. Finally, we prove Theorem \ref{th:lambda-connections via foliations} in 
Section \ref{sc:foliations}. 

\section{Reduction of Theorem \ref{th:lambda-connections} to formal disc}
\label{sc:reduction}

\subsection{Refinement of Theorem \ref{th:lambda-connections}}
Let us now make Theorem \ref{th:lambda-connections} more precise. As before, $X$ is a smooth
curve, $\tX\subset T^*X$ is a smooth spectral curve, $p_\tX:\tX\to X$ is the projection, 
$\{\tx_1,\dots,\tx_n\}\subset\tX$ is the ramification locus of $p_\tX$. Set $\tX_u\eqd \tX-\{\tx_1,\dots,\tx_n\}$,
$X_u\eqd p_\tX(\tX_u)\subset X$. We then have the following commutative diagram of groupoids:
\begin{equation}
\begin{array}{ccc}
\Conn{\tX}&\to&\Conn{\tX_u}\cr
\downarrow&&\downarrow\cr
\Higgs{\tX}&\to&\Higgs{\tX_u},
\end{array}
\label{dg:Conn}
\end{equation}
where $\Higgs{\tX}$ (resp. $\Higgs{\tX_u}$) is the groupoid of Higgs bundles on $X$ (resp. $X_u$)
whose spectral curve is $\tX$ (resp. $\tX_u$), 
and $\Conn{\tX}$ (resp. $\Conn{\tX_u}$) 
is the groupoid of $\C[[\lambda]]$-families of $\lambda$-connections on $X$ (resp. $X_u$)
from Theorem \ref{th:lambda-connections}. In the diagram \eqref{dg:Conn}, the horizontal arrows are functors of 
restriction from $X$ to $X_u$ and the vertical arrows are functors of reduction modulo $\lambda$.

Similarly, the groupoid $\tConn{\tX}$ fits into the commutative diagram
\begin{equation}
\begin{array}{ccc}
\tConn{\tX}&\to&\tConn{\tX_u}\cr
\downarrow&&\downarrow\cr
\tHiggs{\tX}&\to&\tHiggs{\tX_u}.
\end{array}
\label{dg:tConn}
\end{equation}
Here $\tHiggs{\tX}$ (resp. $\tHiggs{\tX_u}$) is the groupoid of line bundles on $\tX$ (resp. $\tX_u$). 

Theorem \ref{th:Hitchin} provides an equivalence $\Higgs{\tX}\iso\tHiggs{\tX}$ and an equivalence $\Higgs{\tX_u}\iso\tHiggs{\tX_u}$).
So we see that the bottom rows of diagrams \eqref{dg:Conn}, \eqref{dg:tConn} are naturally equivalent. 
By the following statement, their upper right
corners are also equivalent:

\begin{theorem} Suppose $X$ is a smooth curve and $p_\tX:\tX\to X$ is an unramified spectral curve. 
The functor $\tConn{\tX}\to\Conn{\tX}$ that sends
$(l,\delta)$ to $((p_\tX)_*(l),(p_\tX)_*(\delta))$ is an equivalence of groupoids. \qed
\label{th:Wasow}
\end{theorem}

This theorem is a bit generalized version of \cite[Proposition 1.2]{DS} (see also \cite[Theorem 25.2]{Wa}) and
can be proved by the same method.
We are now ready to refine the statement of Theorem \ref{th:lambda-connections}:

\begin{theorem}
Let $X$ be a smooth curve and $\tX\subset T^*X$ a smooth spectral curve. 
Consider the fibered products of groupoids:
\begin{align}
\CC{\tX}&\eqd \Conn{\tX_u}\times_{\Higgs{\tX_u}}\Higgs{\tX}
\label{df:CC}\\
\tCC{\tX}&\eqd \tConn{\tX_u}\times_{\tHiggs{\tX_{u}}}\tHiggs{\tX}.
\label{df:tCC}
\end{align}
Notice that Theorems \ref{th:Hitchin} and \ref{th:Wasow} give an equivalence
$\CC{\tX}\iso\tCC{\tX}$.

\begin{enumerate}
\item The functor $\Conn{\tX}\to\CC{\tX}$ induced by \eqref{dg:Conn} 
is fully faithful (so that $\Conn{\tX}$ is a full
subcategory of $\CC{\tX}$).

\item The functor $\tConn{\tX}\to\tCC{\tX}$ induced by \eqref{dg:tConn}
is fully faithful.

\item For a groupoid $\cG$, let $[\cG]$ be the set of isomorphism classes of 
objects of $\cG$. The equivalence $\CC{\tX}\iso\tCC{\tX}$ induces an isomorphism $[\CC{\tX}]\iso[\tCC{\tX}]$.
We claim that this isomorphism identifies the sets $[\Conn{\tX}]\subset[\CC{\tX}]$ and 
$[\tConn{\tX}]\subset[\tCC{\tX}]$. 
\end{enumerate}
\label{th:lambda-connections2}
\end{theorem}

\begin{remark} It is obvious that the restriction functors $\Conn{\tX}\to\Conn{\tX_u}$ and
$\tConn{\tX}\to\tConn{\tX_u}$ are faithful. Therefore, the functors $\Conn{\tX}\to\CC{\tX}$ and
$\tConn{\tX}\to\tCC{\tX}$ are also automatically faithful. 
\end{remark}

Theorem \ref{th:lambda-connections2} claims that the equivalence
$\CC{\tX}\iso\tCC{\tX}$ induces an equivalence $\Conn{\tX}\iso\tConn{\tX}$ that is unique up to
a canonical isomorphism. Therefore, Theorem \ref{th:lambda-connections2} implies Theorem \ref{th:lambda-connections}.

\subsection{}
It is easy to see that all of the above definitions ($\lambda$-connections, Higgs bundles, spectral curves, etc.)
still make sense if $X$ is a formal disc rather then a smooth curve (see Section \ref{sc:formal groupoids} for examples). 
Therefore, we can formulate a `formal' version of Theorem \ref{th:lambda-connections2}:

\begin{theorem} 
Let $X\simeq\spf\C[[z]]$ be a formal disc and $\tX\subset T^*X$ a smooth spectral curve. 
Define $\CC{\tX}$ and $\tCC{\tX}$ by \eqref{df:CC}, \eqref{df:tCC}.
\begin{enumerate}
\item The natural functor $\Conn{\tX}\to\CC{\tX}$ is fully faithful.
\item The natural functor $\tConn{\tX}\to\tCC{\tX}$ is fully faithful.
\item Let us identify $[\CC{\tX}]$ and $[\tCC{\tX}]$ using the equivalence 
$\CC{\tX}\iso\tCC{\tX}$.
We claim that under this identification, $[\Conn{\tX}]\subset[\CC{\tX}]$ corresponds to 
$[\tConn{\tX}]\subset[\tCC{\tX}]$.
\end{enumerate}
\label{th:lambda-connectionsf}
\end{theorem}

We will prove this theorem in the next section. Let us show now that Theorem \ref{th:lambda-connectionsf} implies
Theorem \ref{th:lambda-connections2} (and so also Theorem \ref{th:lambda-connections}).

\begin{proof}[Proof of Theorem \ref{th:lambda-connections2}] Let $X$ be a smooth curve over $\C$, $\tX\subset T^*X$ a smooth 
spectral curve. To simplify the notation, we will assume that $p_\tX:\tX\to X$ is ramified at a single point, $\tx\in\tX$.
Denote by $\tX^\wedge$ the formal completion of $\tX$ at $\tx$ and by $X^\wedge$ the formal completion of $X$ at $x=p_\tX(\tx)$.
Clearly, $X^\wedge$ is a formal disc and $\tX^\wedge$ is a (smooth ramified) spectral curve over $X^\wedge$.

It is a standard fact that the natural diagram 
\begin{equation*}
\begin{array}{ccc}
\Higgs{\tX}&\to&\Higgs{\tX_u}\cr
\downarrow&&\downarrow\cr
\Higgs{\tX^\wedge}&\to&\Higgs{\tX_u^\wedge}
\end{array}
\end{equation*}
is Cartesian; essentially, the claim is that a Higgs bundle on $\tX$ can be glued from a Higgs bundle on $\tX_u$,
a Higgs bundle on $\tX^\wedge$, and an identification of their restrictions to the
punctured disc $\tX_u^\wedge\eqd\tX_u\cap\tX^\wedge$.
The same statement holds for groupoids $\tHiggs{\al}$, $\Conn{\al}$, and $\tConn{\al}$. Now 
Theorem \ref{th:lambda-connections2} follows from Theorem \ref{th:lambda-connectionsf} by diagram chasing: it
suffices to check that the diagrams 
\begin{equation*}
\begin{array}{ccccccc}
\Conn{\tX}&\to&\CC{\tX}&&\tConn{\tX}&\to&\tCC{\tX}\cr
\downarrow&&\downarrow&\quad\text{and}\quad&\downarrow&&\downarrow\cr
\Conn{\tX^\wedge}&\to&\CC{\tX^\wedge}&&\tConn{\tX^\wedge}&\to&\tCC{\tX^\wedge}
\end{array}
\end{equation*}
are Cartesian.
\end{proof}

\subsection{}
\label{sc:formal groupoids}
Let us now describe the groupoids from Theorem \ref{th:lambda-connectionsf} explicitly. 
Set $X=\spf\C[[z]]$. 
The cotangent bundle to $X$ equals $\spf\C[\xi][[z]]$, where $\xi$ is the
vector field $\frac{d}{dz}$ on $X$. A spectral curve $\tX$ over $X$ is given by one equation 
\begin{equation}
\xi^2-t(z)\xi+d(z)=0,\quad (t(z),d(z)\in\C[[z]]). 
\label{eq:formal spectral curve}
\end{equation}
We will only consider the case when $\tX$ is ramified over $X$, because only this case is needed for 
Theorem \ref{th:lambda-connections2} (besides, the unramified case is simply a `formal' version of 
Theorem \ref{th:Wasow}). Since we also want the spectral curve to be smooth, 
$t(z)$ and $d(z)$ must satisfy the following condition:
\begin{equation}
\text{The discriminant }t^2(z)-4d(z)\in\C[[z]]\text{ has a simple zero at }z=0.
\label{cn:non-degeneracy}
\end{equation} 

\begin{notation*}
We denote by $\Omega_X$ the $\C[[z]]$-module of (continuous) differentials of $\C[[z]]$; it is a free $\C[[z]]$-module
generated by $dz$. To simplify the notation, we will write $Ldz$ instead of $L\otimes_{\C[[z]]}\Omega_X$ 
for a $\C[[z]]$-module $L$.  
\end{notation*}

$\Higgs{\tX}$ is the groupoid of pairs $(L,\nabla)$, where $L$ is a rank $2$ free $\C[[z]]$-module and
$\nabla:L\to Ldz$ is a $\C[[z]]$-linear map such that $\tr\nabla=t(z)dz$, $\det(\nabla)=d(z)(dz)^2$.
Here $\Omega_X$ is the free $\C[[z]]$-module with generator $dz$. The groupoid $\Higgs{\tX_u}$ is similar, except
$L$ is a two-dimensional vector space over $\C((z))$.

$\Conn{\tX}$ is the groupoid of pairs $(L,\nabla)$, where $L$ is a rank $2$ free $\C[[z,\lambda]]$-module
and $\nabla:L\to Ldz$ is a ($\C[[\lambda]]$-linear) $\lambda$-connection such that the map 
$\nabla_0:L/\lambda L\to(L/\lambda L)dz$ induced by $\nabla$ satisfies $\tr\nabla_0=t(z)dz$, $\det(\nabla_0)=d(z)(dz)^2$.
The groupoid $\Higgs{\tX_u}$ is similar, except $L$ is a rank $2$ free module over $\C((z))[[\lambda]]$.

$\CC{\tX}$ is the groupoid of triples 
$(L_u,\nabla,L_0)$, where $(L_u,\nabla)\in\Conn{\tX_u}$, and $L_0$ is a $\C[[z]]$-lattice in the vector 
space $L_u/\lambda L_u$ such that $\nabla_0(L_0)\subset L_0dz$. Here the Higgs field 
$\nabla_0:L_u/\lambda L_u\to(L_u/\lambda L_u)dz$ is induced by $\nabla$. 

$\tHiggs{\tX}$ (resp. $\tHiggs{\tX_u}$) is the groupoid of rank $1$ free $\C[[\tz]]$-modules (resp. one-dimensional
$\C((\tz))$-vector spaces), where $\tz$ is a formal coordinate on 
$\tX\simeq\spf\C[[\tz]]$. For instance, we can set 
\begin{equation}
\tz=\sqrt{t(z)^2-4d(z)}.
\label{eq:tz}
\end{equation} 

$\tConn{\tX}$ is the groupoid of pairs $(l,\delta)$, where $l$ is a rank $1$ free $\C[[\tz,\lambda]]$-module, and
$\delta:l\to \tz^{-1}l d\tz$ is a ($\C[[\lambda]]$-linear) $\lambda$-connection such that 
$\res\delta:l/\tz l\to(\tz^{-1}l/l)d\tz$ equals $-\lambda/2$, and
the induced map $\delta_0:l/\lambda l\to \tz^{-1}(l/\lambda l)d\tz$ equals the natural $1$-form $\mu\in\Omega_\tX$. 
Notice that if $\tz$ is given by \eqref{eq:tz}, then $\mu=(-t(z)+\tz)dz/2$. Similarly, $\tConn{\tX_u}$ is
the groupoid of pairs $(l,\delta)$, where $l$ is a rank $1$ free $\C((\tz))[[\lambda]]$-module,
and $\delta:l\to ld\tz$ is a $\lambda$-connection such that the induced map $\delta_0:l/\lambda l\to(l/\lambda l)d\tz$
equals $\mu$.

Finally, $\tCC{\tX}$ is the groupoid of triples $(l_u,\delta,l_0)$, where $(l_u,\delta)\in\tConn{\tX_u}$, 
and $l_0$ is a $\C[[z]]$-lattice in the vector space $l_u/\lambda l_u$.

Let us now describe the natural functors between these groupoids. The functor $\Conn{\tX}\to\CC{\tX}$
(induced by the diagram \eqref{df:CC})
sends $(L,\nabla)\in\Conn{\tX}$ to $(L\otimes_{\C[[z]]}\C((z)),\nabla,L/\lambda L)\in\CC{\tX}$. Similarly,
the functor $\tConn{\tX}\to\tCC{\tX}$ (induced by the diagram \eqref{df:tCC})
sends $(l,\delta)\in\tConn{\tX}$ to $(l\otimes_{\C[[z]]}\C((z)),\delta,l/\lambda l)\in\tCC{\tX}$.

Finally, let us describe the equivalence $\tCC{\tX}\iso\CC{\tX}$. Given $(l_u,\delta,l_0)\in\tCC{\tX}$, we can consider
$l_u$ as a $\C((z))[[\lambda]]$-module using the embedding $\C[[z]]\hookrightarrow\C[[\tz]]$. Similarly, $l_0$ can 
be viewed as a $\C[[z]]$-lattice in the two-dimensional $\C((z))$-space $l_u/\lambda l_u$. We have
$\Omega_X\otimes\C((\tz))=\Omega_\tX\otimes\C((\tz))$, therefore $l_udz=l_ud\tz$ and we can view $\delta$ as
a $\lambda$-connection on the $\C((z))[[\lambda]]$-module $l_u$. In this sense, the functor
$\tCC{\tX}\iso\CC{\tX}$ is forgetful: it sends $(l_u,\delta,l_0)$ to the same triple $(l_u,\delta,l_0)$, but considered
over $\C[[z]]$.

We will also need the following property of the equivalence $\tCC{\tX}\iso\CC{\tX}$. Take any $(l_u,\delta,l_0)\in\tCC{\tX}$,
and let $(L_u,\nabla,L_0)$ be its image in $\CC{\tX}$ (so that $l_u$ and $L_u$ are identified as $\C[[z]]$-modules).
Set $\tL_u\eqd L_u\otimes\C((\tz))$. We then obtain a natural isomorphism
\begin{equation*}
\phi:\tL_u\iso l_u\oplus\sigma^* l_u,
\end{equation*}
where $\sigma$ is the non-trivial element of the Galois 
group $\gal(\C((\tz))/\C((z)))$. Notice that $\nabla$ induces a $\lambda$-connection $\tnabla:\tL_u\to\tL_u d\tz$,
and $\delta$ induces a $\lambda$-connection $\delta\oplus\sigma^*\delta$ on $l_u\oplus\sigma^* l_u$. 
Besides, $\tL_0\eqd L_0\otimes\C((\tz))$ is a $\C[[\tz]]$-lattice in $\tL_u/\lambda\tL_u$,
and $l_0\oplus\sigma^* l_0$ is a $\C[[\tz]]$-lattice in $(l_u\oplus\sigma^* l_u)/\lambda(l_u\oplus\sigma^* l_u)$. The
following lemma is easy to prove:

\begin{lemma} 
\begin{enumerate}
\item $\phi$ agrees with the $\lambda$-connections: $\phi\circ\tnabla=(\delta\oplus\sigma^*\delta)\circ\phi$.
\item The isomorphism 
\begin{equation*}
\phi_0:\tL_u/\lambda\tL_u\iso (l_u\oplus\sigma^* l_u)/\lambda(l_u\oplus\sigma^* l_u)
\end{equation*} 
(induced by $\phi$) satisfies $\phi_0(\tL_0)\subset l_0\oplus\sigma^* l_0$.
\item The image $\phi_0(\tL_0)\subset l_0\oplus\sigma^* l_0$ is the kernel of the composition
\begin{equation*}
l_0\oplus\sigma^* l_0\to l_0/\tz l_0\oplus(\sigma^* l_0)/(\tz\sigma^* l_0)=(l_0/\tz l_0)^2\to(l_0/\tz l_0),
\end{equation*}
where the last map is $(x,y)\mapsto x+y$. In particular, the quotient $(l_0\oplus\sigma^* l_0)/\phi_0(\tL_0)$
is one-dimensional.
\end{enumerate}
\label{lm:equi} \qed
\end{lemma}

\section{$\lambda$-connections on formal disc: Theorem \ref{th:lambda-connectionsf}}
\label{sc:proof}

Let us keep the notation of Section \ref{sc:formal groupoids}. So $X=\spf\C[[z]]$, $\tX$ is given by \eqref{eq:formal spectral curve},
and $t(z),d(z)\in\C[[z]]$ satisfy \eqref{cn:non-degeneracy}.

\subsection{Proof of Theorem \ref{th:lambda-connectionsf}(1)}
\label{sc:proof1}

\newcommand\grc{G}		
\newcommand\grf[1]{G_{#1}}	
\newcommand\gauge[2]{Gauge_\lambda(#1,#2)}

\begin{lemma}
Suppose $A_0(z)\in\gl_2(\C[[z]])$ is such that $t(z)=\tr A_0(z)$, $d(z)=\det A_0(z)$.
\begin{enumerate}
\item $A_0(z)$ is a regular element of $\gl_2(\C[[z]])$.
\item Suppose $B(z)\in\gl_2(\C((z)))$ satisfies 
$[B(z),A_0(z)]=0$ and 
\begin{equation*}
\frac{dB(z)}{dz}\in\gl_2(\C[[z]])+[A_0(z),\gl_2(\C((z)))].
\end{equation*}
Then $B(z)\in\gl_2(\C[[z]])$.
\end{enumerate}
\label{lm:non-degeneracy}
\end{lemma}

\newcommand\pz[1]{\frac{d#1}{dz}}
\begin{proof}
The first statement of the lemma is almost obvious; let us prove the second one. 
Notice that $\tr\pz{B},\tr\left(\pz{B}A_0\right)\in\C[[z]]$. 
Since $[B(z),A_0(z)]=0$, and $A_0(z)$ is regular, we can write 
\begin{equation*}
B(z)=f(z)+g(z)A_0(z)
\end{equation*} 
for some $f,g\in\C((z))$. We have
\begin{align}
\tr\pz{B}&=2\pz{f}+\pz{g}\tr A_0+g\tr\pz{A_0}=2\pz{f}+\pz{g}t+g\pz{t},
\label{eq:tr1}\\
\intertext{and}
\begin{split}
\tr\left(\pz{B}A_0\right)
&=\pz{f}\tr A_0+\pz{g}\tr(A_0^2)+g\tr\left(\pz{A_0}A_0\right)\\
&=\pz{f}t+\pz{g}(t^2-2d)+g\left(t\pz{t}-d\right).
\end{split}\\
\intertext{Therefore,}
\tr\left(\pz{B}(2A_0-t)\right)&=\pz{g}(t^2-4d)+\frac{1}{2}g\pz{(t^2-4d)}.
\label{eq:tr2}
\end{align}  
Clearly, \eqref{eq:tr2} belongs to $\C[[z]]$; by looking at the leading term of the Laurent expansion of \eqref{eq:tr2},
it is easy to see that $g\in\C[[z]]$. Now \eqref{eq:tr1} implies $\pz{f}\in\C[[z]]$, so $f\in\C[[z]]$. Finally,
$B=f+gA_0\in\gl_2(\C[[z]])$.
\end{proof}

\begin{remark} Let us reformulate Lemma \ref{lm:non-degeneracy}(2). Let $\cnt\eqd\ker(\ad A_0)$ and 
$\ccnt\eqd\coker(\ad A_0)$ be the centralizer and the `co-centralizer' of $A_0\in\gl_2(\C[[z]])$, respectively.
Then $\cnt$ and $\ccnt$ are rank $2$ free $\C[[z]]$-modules, and $\cnt\otimes\C((z))$ and $\ccnt\otimes\C((z))$ are 
the centralizer and the co-centralizer of $A_0$ in $\gl_2(\C((z)))$, respectively. Let $D$ be the composition
\begin{equation*}
\cnt\otimes\C((z))\hookrightarrow\gl_2(\C((z)))\xrightarrow{d}\gl_2(\C((z)))dz\to\ccnt\otimes\C((z))dz.
\label{eq:Dcnt}
\end{equation*}    
Then Lemma \ref{lm:non-degeneracy}(2) claims that $D^{-1}(\ccnt dz)=\cnt$.

Note also that the composition $\cnt\to\gl_2(\C[[z]])\to\ccnt$ identifies $\cnt$ with a submodule of $\ccnt$. 
In this manner, $D$ can be viewed as a connection (with a pole at $z=0$) on the $\C[[z]]$-module $\cnt$. 
The monodromy of $D$ has eigenvalues $1$ and $-1$; this follows from Theorem 11.6 of \cite{DG} (actually,
Step 2 in the proof of \cite[Proposition 12.5]{DG} suffices). It is not hard to derive Lemma \ref{lm:non-degeneracy}(2)
from this observation.
\end{remark}

\begin{lemma}
\label{lm:lambda-gauge}
Suppose $A(z,\lambda)\in\gl_2(\C[[z,\lambda]])$ is such that $\tr(A(z,0))=t(z)$, $\det(A(z,0))=d(z)$, and
$R(z,\lambda)\in\GL_2(\C((z))[[\lambda]])$ satisfies $R(z,0)\in\GL_2(\C[[z]])$. Set
\begin{equation}
\gauge{A}{R}
\eqd R^{-1}AR+\lambda R^{-1}\frac{dR}{dz}
\in\gl_2(\C((z))[[\lambda]]),
\label{eq:lambda-gauge}
\end{equation}
and suppose $\gauge{A}{R}\in\gl_2(\C[[z,\lambda]])$. Then $R(z,\lambda)\in\GL_2(\C[[z,\lambda]])$.
\end{lemma}

\begin{remark}
$\gauge{A}{R}$ is a $\lambda$-version of gauge transform: we can rewrite \eqref{eq:lambda-gauge} 
as
\begin{equation*}
\left(\lambda\frac{d}{dz}+\gauge{A}{R}\right)=R^{-1}(z,\lambda)\left(\lambda\frac{d}{dz}+A(z,\lambda)\right) R(z,\lambda).
\end{equation*}
\end{remark}

\begin{proof}
Assume the converse. Let us expand 
\begin{equation*}
A(z,\lambda)=\sum_{i\ge0} A_i(z)\lambda^i,\quad
\gauge{A}{R}=\sum_{i\ge0} A'_i(z)\lambda^i,\quad
R(z,\lambda)=\sum_{i\ge0} R_i(z)\lambda^i.
\end{equation*}
Note that $A_i(z),A'_i(z)\in\gl_2(\C[[z]])$ for $i\ge0$, $R_0(z)\in\GL_2(\C[[z]])$, and $R_i(z)\in\gl_2(\C((z)))$ for $i>0$.

Let $j>0$ be the minimal index such that $R_j(z)\not\in\gl_2(\C[[z]])$. Set 
\begin{equation*}
R'\eqd\sum_{i=0}^{j-1}R_j(z)\lambda^j\in\GL_2(\C[[z,\lambda]]),\quad B\eqd\gauge{A}{R'},\quad Q\eqd(R')^{-1}R.
\end{equation*} 
Then $\gauge{A}{R}=\gauge{B}{Q}$. Replacing $A$ with $B$ and
$R$ with $Q$, we can assume without loss of generality that $R'=1$, that is, $R_0(z)=1$ and $R_i(z)=0$ 
for $0<i<j$. 

Taking the coefficient of $\lambda^j$ in \eqref{eq:lambda-gauge}, we now obtain
\begin{equation*}
A'_j(z)=A_j(z)+[A_0(z),R_j(z)].
\end{equation*}
Therefore, $[A_0(z),R_j(z)]\in\gl_2(\C[[z]])$. By Lemma \ref{lm:non-degeneracy}(1), 
we know that $A_0(z)\in\gl_2(\C[[z]])$ is a
regular element, and so we can write $R_j(z)=S(z)+T(z)$ for some $S(z)\in\gl_2(\C[[z]])$,
$T(z)\in\gl_2(\C((z)))$ such that $[T(z),A_0(z)]=0$. Replacing $A(z,\lambda)$ with $\gauge{A}{1+S(z)\lambda^j}$
and $R(z,\lambda)$ with $(1+S(z)\lambda^j)^{-1}R$, we can assume $S(z)=0$, so that $[A_0(z),R_j(z)]=0$.

Now taking the coefficient of $\lambda^{j+1}$ in \eqref{eq:lambda-gauge}, we obtain
\begin{equation*}
A'_{j+1}(z)=A_{j+1}(z)+\frac{dR_j}{dz}(z)+[A_0(z),R_{j+1}(z)]+[A_1(z),R_j(z)]-[R_1(z)A_0(z),R_j(z)]
\end{equation*}
(the last term is non-zero only if $j=1$). Recall that $[A_0(z),R_j(z)]=0$ and $A_0(z)$ is a regular element;
therefore, 
\begin{equation*}
[A_1(z),R_j(z)]-[R_1(z)A_0(z),R_j(z)]\in[A_0(z),\gl_2(\C((z)))].
\end{equation*}
Lemma \ref{lm:non-degeneracy}(2) implies $R_j(z)\in\gl_2(\C[[z]])$, which contradicts our assumption.
\end{proof}

Lemma \ref{lm:lambda-gauge} can be reformulated in terms of groupoid $\CC{\tX}$:

\begin{corollary} For any $(L_u,\nabla,L_0)\in\CC{\tX}$, there exists at most one $\C[[z,\lambda]]$-lattice
$L\subset L_u$ such that $L_0=L/\lambda L\subset L_u/\lambda L_u$ and $\nabla(L)\subset Ldz$.
\label{co:lambda-gauge}
\end{corollary}
\begin{proof}
Suppose $L\subset L_u$ is such a lattice. Choose a trivialization 
$\iota:(\C[[z,\lambda]])^2\iso L$, and define the \emph{connection matrix} of $\nabla$ to be 
$A\in\gl_2(\C[[z,\lambda]])$ such that
\begin{equation*}
\left(\lambda\frac{d}{dz}+A(z,\lambda)\right)dz=\iota^{-1}\nabla\iota.
\end{equation*}

Now let $L'\subset L_u$ be another such lattice; denote by $A'\in\gl_2(\C[[z,\lambda]]$ the connection matrix
corresponding to some trivialization $\iota':(\C[[z,\lambda]])^2\iso L'$.
Clearly, $A'$ is obtained from $A$ by the `$\lambda$-gauge transform' \eqref{eq:lambda-gauge}:
$A'=\gauge{A}{R}$ for
$R(z,\lambda)\eqd \iota^{-1}\iota'\in\GL_2(\C((z))[[\lambda]])$. As 
$L/\lambda L=L'/\lambda L'\subset L_u/\lambda L_u$, we see that $R(z,0)\in\GL_2(\C[[z]])$, so that
the assumptions of Lemma \ref{lm:lambda-gauge} are satisfied. Therefore, $R(z,\lambda)\in\GL_2(\C[[z,\lambda]])$,
and $L=L'$.
\end{proof}

Recall that for a functor $\cF:\cG_1\to\cG_2$ between groupoids, the \emph{(essential) fiber} of $\cF$
over $\gamma_2\in\cG_2$ is the groupoid of pairs
\begin{equation*}
\{(\gamma_1,f)\st\gamma_1\in\cG_1,f:\cF(\gamma_1)\iso\gamma_2\}.
\end{equation*} 
Given $(L_u,\nabla,L_0)\in\CC{\tX}$, the set of all $\C[[z,\lambda]]$-lattices $L$ as in Corollary \ref{co:lambda-gauge}
is equivalent to the fiber of the functor $\Conn{\tX}\to\CC{\tX}$; 
Corollary \ref{co:lambda-gauge} claims the fiber is either empty or equivalent to a one-element set. Now Theorem
\ref{th:lambda-connectionsf}(1) is implied by the following simple lemma:

\begin{lemma} Let $\cF:\cG_1\to\cG_2$ be a functor between groupoids. 
\begin{enumerate}
\item $\cF$ is faithful if and only if the fiber of $\cF$ over any object of $\cG_2$ is discrete.

\item $\cF$ is fully faithful if and only if the fiber of $\cF$ over any object of $\cG_2$ is either empty
or equivalent to a one-element set.
\end{enumerate}
\qed
\label{lm:groupoids}
\end{lemma}

\subsection{Proof of Theorem \ref{th:lambda-connectionsf}(2)} \label{sc:proof2}

The second statement of Theorem \ref{th:lambda-connectionsf} is proved similarly to its first statement.
Actually, Theorem \ref{th:lambda-connectionsf}(2) is simpler, because it deals with `abelian' objects 
(line bundles); for instance,
the $\lambda$-gauge transformation \eqref{eq:lambda-gauge} simplifies. 
Finally, notice that $\lambda$-connections on line bundles can be 
reduced to
ordinary connections (as in Remark \ref{rm:lambda-connections on line bundles}), so all results of this section are 
more or less classical. 

For any $(l_u,\delta,l_0)\in\tCC{\tX}$, choose a trivialization $\iota:\C((\tz))[[\lambda]]\iso l_u$
such that the induced map $\C((\tz))\iso l_u/\lambda l_u$ identifies $\C[[\tz]]$ and $l_0$. We will say that
$\iota$ respects $l_0$; clearly, such $\iota$
always exists. Denote by $a(\tz,\lambda)\in\C((\tz))[[\lambda]]$ the connection matrix of 
$\delta$ with respect to $\iota$:
\begin{equation*}
\left(\lambda\frac{d}{d\tz}+a(\tz,\lambda)\right)d\tz=\iota^{-1}\delta\iota.
\end{equation*}
Recall now that the map $\delta_0:l/\lambda l\to(l/\lambda l)d\tz$ induced by $\delta$ equals $\mu$ 
(where $\mu\in\Omega_\tX$ is the canonical 1-form on $\tX$). 
Therefore, $a(\tz,0)=\mu(d\tz)^{-1}$. 

\begin{lemma} Let $(l_u,\delta,l_0)$, $\iota$, and $a(\tz,\lambda)$ be as above. Denote by $Y$ 
the set of formal series $r(\tz,\lambda)\in\C((z))[[\lambda]]$ that satisfy the following conditions:
\begin{gather}
r(\tz,0)\in\C[[\tz]]^\times\label{cn:r1}\\
\text{Set }a'\eqd a+\lambda r^{-1}\frac{dr}{d\tz}.\text{ Then }
a'\in\tz^{-1}\C[[\tz,\lambda]],\text{ and }\res_{\tz=0}a'=-\lambda/2.\label{cn:r2}
\end{gather}
Let $\C[[z,\lambda]]^\times$ be the group of invertible Taylor power series of two variables; we let
$\C[[z,\lambda]]^\times$ act on $Y$ by multiplication. Then the fiber of 
the functor $\tConn{\tX}\to\tCC{\tX}$ over $(l_u,\delta,l_0)$ is equivalent 
to the quotient set $Y/\C[[z,\lambda]]^\times$.
\label{lm:fiber}
\end{lemma}
\begin{proof}
The fiber of the functor is equivalent to the set of all $\C[[\tz,\lambda]]$-lattices $l\subset l_u$ such that
\begin{gather}
l/\lambda l=l_0\subset l_u/\lambda l_u\label{cn:l1}\\
\intertext{and} 
(l,\delta)\in\tConn{\tX}.\label{cn:l2}
\end{gather} 
Any such lattice $l$ can be written as 
\begin{equation*}
l=r\iota(\C[[\tz,\lambda]])
\end{equation*}
for some $r(\tz,\lambda)\in(\C((\tz))[[\lambda]])^\times$. Notice that $r$ is unique up to multiplication by
an invertible element of $\C[[\tz,\lambda]]$. Also, $a'(\tz,\lambda)=\gauge{a}{r}$ is the connection matrix 
of $\delta$ with respect to the trivialization $r\iota:\C[[\tz,\lambda]]\iso l$. Now it is clear that
\eqref{cn:r1} and \eqref{cn:r2} are equivalent to \eqref{cn:l1} and \eqref{cn:l2}, respectively. This completes 
the proof. 
\end{proof}

\begin{lemma} Suppose $a(\tz,\lambda)\in\C((\tz))[[\lambda]]$ satisfies $a(\tz,0)=\mu(d\tz)^{-1}$. Define 
the set $Y$ as in Lemma \ref{lm:fiber}. Then either $Y$ is empty or it consists of a single 
$\C[[\tz,\lambda]]^\times$-orbit. $Y$ is non-empty if and only if $a(\tz,\lambda)$ satisfies the following
condition:
\begin{equation}
a_1(\tz)\in\tz^{-1}\C[[\tz]]\text{ and }\res_{\tz=0}a(\tz,\lambda)=-\lambda/2.
\label{cn:non-empty}
\end{equation} 
Here $a_1(\tz)\eqd\frac{\partial a(\tz,\lambda)}{\partial\lambda}$ is the coefficient of $\lambda$ in 
$a(\tz,\lambda)$.
\end{lemma}

\begin{proof} Denote by $Z$ the set of all series $r(\tz,\lambda)\in\C((\tz))[[\lambda]]$ that satisfy \eqref{cn:r1},
and consider the map
\begin{equation*}
d\log:Z\to\C((\tz))[[\lambda]]:r(\tz,\lambda)\mapsto r(\tz,\lambda)^{-1}\frac{dr(\tz,\lambda)}{d\tz}.
\end{equation*}
Notice that $Z$ is a group under multiplication and $d\log$ is a group homomorphism. To prove the lemma, we need to
verify two properties of $d\log$:
\begin{gather*}
d\log(Z)=\{s(\tz,\lambda)\in\C((\tz))[[\lambda]]\st s(\tz,0)\in\C[[\tz]],\res_{\tz=0}s(\tz,0)=0\}\\
(d\log)^{-1}(\C[[\tz,\lambda]])=\C[[\tz,\lambda]]\subset Z.
\end{gather*}
Both properties are almost obvious, especially if one notices that for $r(\tz,\lambda)\in Z$, the expression 
$\log(r)\in\C((\tz))[[\lambda]]$ makes sense.
\end{proof}

\begin{corollary} Let $(l_u,\delta,l_0)$, $\iota$, and $a(\tz,\lambda)$ be as in Lemma \ref{lm:fiber}, and let 
$a_1(\tz)\eqd\frac{\partial a(\tz,\lambda)}{\partial\lambda}$ be the coefficient of $\lambda$ in $a(\tz,\lambda)$.
Then the fiber of the functor $\tConn{\tX}\to\tCC{\tX}$ over $(l_u,\delta,l_0)$ is either 
equivalent to a one-point set or empty. The fiber is non-empty if and only if \eqref{cn:non-empty} holds. \qed
\label{co:fibertilde}
\end{corollary}

By Lemma \ref{lm:groupoids}, 
Corollary \ref{co:fibertilde} is equivalent to Theorem \ref{th:lambda-connectionsf}(2). Let us 
also prove the following lemma, which is used in the next section.

\begin{lemma} All objects of $\tConn{\tX}$ are isomorphic. 
\label{lm:tConn}
\end{lemma}
\begin{proof}
Take any $(l,\delta), (l',\delta')\in\Conn{\tX}$, and let us choose trivializations $\iota:\C[[\tz,\lambda]]\iso l$,
$\iota':\C[[\tz,\lambda]]\iso l'$. Denote by $a(\tz,\lambda),a'(\tz,\lambda)$ the connection matrices of $\delta$ and
$\delta'$, respectively. Notice that $a-a'\in\lambda\C[[\tz,\lambda]]$, because $\res_{\tz=0}a=\res_{\tz=0}a'$ and
$a(\tz,0)=a'(\tz,0)$. Therefore, there exists $r\in(\C[[\tz,\lambda]])^\times$ such that 
\begin{equation*}
a'(\tz,\lambda)-a(\tz,\lambda)=\lambda r^{-1}(\tz,\lambda)\frac{d r(\tz,\lambda)}{d\tz}=\lambda d\log(r(\tz,\lambda)).
\end{equation*}
The map $r(\tz,\lambda)\iota'\iota^{-1}$ provides an isomorphism $(l,\delta)\iso(l',\delta')$.
\end{proof}

\subsection{Proof of Theorem \ref{th:lambda-connectionsf}(3)} \label{sc:proof3}

\begin{lemma} Suppose $A(z,\lambda)=\sum_{i\ge 0}A_i(z)\lambda^i\in\gl_2(\C((z))[[\lambda]])$ 
and $R(z,\lambda)\in\GL_2(\C((z))[[\lambda]])$ satisfy the following conditions:
\begin{enumerate}
\item $A(z,\lambda)$ is diagonal;
\item $R(z,0)\in\gl_2(\C[[z]])$;
\item $\det R(z,0)\in\C[[z]]$ has a first-order zero at $z=0$;
\item $\tA:=\gauge{A}{R}\in\gl_2(\C[[z,\lambda]])$ (for the definition of $\gauge{A}{R}$, see \eqref{eq:lambda-gauge}).
\end{enumerate}
Then $A_1(z)\in z^{-1}\gl_2(\C[[z]])$.
\label{lm:regular}
\end{lemma}
\begin{proof}
Set $S(z,\lambda)=R(z,\lambda)R(z,0)^{-1}$, $B(z,\lambda)=\gauge{A}{S}$. Then $B(z,\lambda)=\gauge{\tA}{R(z,0)}$;
this clearly implies $B(z,\lambda)\in z^{-1}\gl_2(\C[[z,\lambda]])$.

Let us expand
\begin{equation*}
B(z,\lambda)=\sum_{i\ge0} B_i(z)\lambda^i,\quad
S(z,\lambda)=\sum_{i\ge0} S_i(z)\lambda^i,
\end{equation*}
where $B_i(z)\in z^{-1}\gl_2(\C[[z]])$, $S_i(z)\in\gl_2(\C((z)))$ for $i\ge 0$, and $S_0(z)=1$. Taking the coefficient
of $\lambda$ in the identity $B=\gauge{A}{S}$, we obtain
\begin{equation*}
B_1(z)=A_1(z)+[A_0(z),S_1(z)].
\end{equation*}
Note that the diagonal entries of $[A_0(z),S_1(z)]$ vanish (because $A_0(z)$ is diagonal), and therefore $A_1(z)$ and
$B_1(z)$ have the same diagonal entries. Therefore, $A_1(z)\in z^{-1}\C[[z]]$ (because $A_1(z)$ is diagonal). 
\end{proof}

\begin{lemma} Let $A$, $R$, and $\tA$ be as in Lemma \ref{lm:regular}. Then $\res_{z=0}\tr A(z,\lambda)=-\lambda$.
\label{lm:residue}
\end{lemma}
\begin{proof}
It is easy to see that
\begin{equation*}
\tr(\tA)=\tr(A)+\lambda\frac{d(\det R(z))}{dz}(\det R(z))^{-1}.
\end{equation*}
Notice that $\det R(z)=z f(z,\lambda)$, where $f(z,\lambda)\in\C((z))[[\lambda]]$ and $f(z,0)\in\C[[z]]^\times$.
Therefore, $\ln f(z,\lambda)$ is well defined, and we can write 
\begin{equation*}
\frac{d(\det R(z))}{dz}(\det R(z))^{-1}=z^{-1}+\frac{df}{dz}f^{-1}=z^{-1}+\frac{d\ln f}{dz}.
\end{equation*}
Hence $\res_{z=0}\tr A(z,\lambda)=\res_{z=0}\tr\tA(z,\lambda)-\lambda=-\lambda$.
\end{proof}

\begin{proposition}
Let $(L_u,\nabla,L_0)\in\CC{\tX}$ be the image of $(L,\nabla)\in\Conn{\tX}$ (so $L_u=L\otimes\C((z))$ and
$L_0=L/\lambda L\subset L_u/\lambda L_u$), and let $(l_u,\delta,l_0)$ be the corresponding object of
$\tCC{\tX}$. Then $(l_u,\delta,l_0)\in\tCC{\tX}$ is isomorphic to the image of an object of
$\tConn{\tX}$ (that is, $(l_u,\delta,l_0)$ belongs to the essential image of $\tConn{\tX}\to\tCC{\tX}$).
\label{pp:image}
\end{proposition}
\begin{proof}
Choose a trivialization $\C((\tz))[[\lambda]]\iso l_u$ that respects $l_0$ (as in Section \ref{sc:proof2}). Denote by
$a(\tz,\lambda)=\sum_{i\ge 0} a_i(\tz)\lambda^i\in\C((\tz))[[\lambda]]$ the matrix of $\delta$ in this trivialization.
According to Corollary \ref{co:fibertilde}, we need to verify \eqref{cn:non-empty} to prove the proposition.
We will do this by using Lemma \ref{lm:equi}.

Let $\sigma$ be the non-trivial element of the Galois group $\gal(\C((\tz))/\C((z)))$.
The trivialization $\C((\tz))[[\lambda]]\iso l_u$ induces a trivialization 
$(\C((\tz))[[\lambda]])^2\iso l_u\oplus\sigma^* l_u$; let 
$A(\tz,\lambda)=\sum_{i\ge 0}A_i(\tz)\lambda^i\in\gl_2(\C((\tz))[[\lambda]])$ be the matrix of the connection 
$\delta\oplus\sigma^*\delta$ with respect to this trivialization. Note that $A(\tz,\lambda)$ is a diagonal matrix,
and one of its entries equals $a(\tz,\lambda)$. Besides, 
\begin{equation*}
\res_{\tz=0}\tr A(\tz,\lambda)=\res_{\tz=0}\delta+\res_{\tz=0}\sigma^*\delta=2\res_{\tz=0}\delta=2\res_{\tz=0}a(\tz,\lambda).
\end{equation*}
Therefore, it suffices to verify that $A_1(\tz)\in\tz^{-1}\gl_2(\C[[\tz]])$ and $\res_{\tz=0}\tr A(\tz,\lambda)=-\lambda$.

Let us choose a trivialization $\iota:(\C[[\tz,\lambda]])^2\iso\tL\eqd L\otimes_{\C[[z]]}\C[[\tz]]$, and let 
$\tA(\tz,\lambda)\in\gl_2(\C[[\tz,\lambda]])$ be the matrix of the connection $\tnabla:\tL\to\tL d\tz$ (induced by $\nabla$)
with respect to $\iota$. From Lemma \ref{lm:equi}, we obtain a morphism $\phi:\tL\to l\oplus l_u$ that respects the 
$\lambda$-connections. Therefore, $\tA=\gauge{A}{R}$, where $R$ is the matrix of $\phi$. Now the required properties
of $A(\tz,\lambda)$ follow from Lemmas \ref{lm:regular} and \ref{lm:residue}.
\end{proof}

Now Theorem \ref{th:lambda-connectionsf}(3) easily follows. Indeed, Proposition \ref{pp:image} implies
that under the isomorphism $[\CC{\tX}]\iso[\tCC{\tX}]$, 
the image of $[\Conn{\tX}]\subset[\CC{\tX}]$ is contained in $[\tConn{\tX}]\subset[\CC{\tX}]$. But
$[\tConn{\tX}]$ is a one-element set (Lemma \ref{lm:tConn}), and $[\Conn{\tX}]$ is obviously not empty. 
Therefore, the isomorphism identifies
the two sets.

\section{Proof of Theorem \ref{th:lambda-connections via foliations}}

\subsection{}
Let us start with a simple observation about foliations on formal schemes:

\begin{definition} A \emph{$\lambda$-adic formal scheme} is a formal scheme $S$ together with a 
function $\lambda\in H^0(S,O_S)$ such that the zero locus of $\lambda^{i+1}$ is a subscheme $S_i\subset S$ and
$S=\limind S_i$. A $\lambda$-adic formal scheme $S$ is \emph{flat} if $S_i$ is flat over $\C[\lambda]/(\lambda^{i+1})$ 
for all $i\ge 0$, or, equivalently, if $\lambda\in O_S$ is not a zero divisor. 
Finally, a $\lambda$-adic formal scheme $S$ is \emph{smooth} if $S_i$ is smooth over $\C[\lambda]/(\lambda^{i+1})$
for all $i\ge0$, or, equivalently, if $S$ is flat and $S_0$ is smooth over $\C$.
\end{definition}

\begin{example}\label{ex:lambda} 
For an arbitrary $\C$-scheme $S$, set $S[[\lambda]]\eqd \limind S\times\spec\C[\lambda]/(\lambda^i)$
(as in Definition \ref{df:C[[lambda]]-family}). Then $S[[\lambda]]$ is a flat $\lambda$-adic formal scheme;
it is smooth if and only if $S$ is smooth.
\end{example}

\begin{lemma} Let $Y$ and $Z$ be smooth $\lambda$-adic formal schemes, and $\Phi:Y\to Z$ a morphism
over $\C[[\lambda]]$ (that is, $\Phi^*(\lambda)=\lambda$). Denote by $Y_0\subset Y$ and $Z_0\subset Z$ the zero loci of $\lambda$.
\begin{enumerate}
\item If the restriction of $\Phi$ to $Y_0$ is smooth, then so is $\Phi$.
\item Suppose $\Phi$ is smooth, and let $\zeta\eqd \ker(d\Phi)\subset T Y$ be the foliation corresponding to the fibration
$\Phi:Y\to Z$. Suppose that the quotient $Y_0/\zeta$ exists and coincides with $Z_0$ (that is, the restriction
$\Phi|_{Y_0}:Y_0\to Z_0$ has connected non-empty fibers). Then the quotient $Y/\zeta$ exists and coincides with $Z$.
\end{enumerate}
\qed
\label{lm:formal foliations}
\end{lemma}

Our proof of Theorem \ref{th:lambda-connections via foliations} is divided into the following steps:
\begin{itemize}
\item Construction of a map $\Phi:\MS[[\lambda]]\to\Cnfp$.
\item Verification that $\Phi$ satisfies the assumptions of Lemma \ref{lm:formal foliations}. 
Therefore, $\Cnfp=\MS[[\lambda]]/\zeta$, where the foliation $\zeta$ equals $\ker(d\Phi)$. 
\item Verification that $\zeta=\zeta_\lambda$. 
\end{itemize}

\subsection{Construction of $\Phi:\MS[[\lambda]]\to\Cnfp$}
\label{sc:foliations}

Even though we formulated Theorems \ref{th:lambda-connections}, \ref{th:lambda-connections2}, and
\ref{th:lambda-connectionsf} for $\C[[\lambda]]$-families
of $\lambda$-connections, the same proof works for $K[[\lambda]]$-families of $\lambda$-connections
on a smooth curve $X/K$, where $K$ is an 
arbitrary $\C$-algebra. The only change is that in Section \ref{sc:formal groupoids}, the modules
over various algebras ($K[[z,\lambda]]$, $K((\tz))$, ans so on) have to be locally free, rather than
free as before; accordingly, all calculations in Sections \ref{sc:proof} have to
be done locally on $K$.

\begin{remark} Of course, the theorems also hold for families parametrized by $S[[\lambda]]$ where
$S$ is a $\C$-scheme (or a stack); indeed, the statements are local on $S$. 
Even more generally, we can consider families parametrized by a flat $\lambda$-adic formal scheme. 
The details of this generalization are left to the reader. Notice 
that once Theorem \ref{th:lambda-connections via foliations} is proved, such statements become 
almost obvious.
\end{remark}

Recall now that $\MS$ is the moduli stack of triples $(\tX,l,\partial)$, where 
$\tX\subset T^*X$ is a smooth spectral curve, $l$ is a line bundle
on $\tX$, and $\partial:l\to l\otimes\Omega(\tx_1+\dots+\tx_n)$ is a connection whose residues 
at $\tx_1,\dots,\tx_n$ (the ramification locus of $p_\tX:\tX\to X$) equal $-1/2$.

Since $\MS$ is a moduli stack, it carries a universal family; let us denote it by $(\tXS,\lS,\partialS)$. Here 
$\tXS\subset (T^*X)\times\MS$ is an $\MS$-family of smooth spectral curves, $\lS$ is a line bundle on $\tXS$, and 
$\partialS:\lS\to\lS\otimes\Omega_{\tXS/\MS}(D_r)$ is a connections with pole at $D_r$, the ramification divisor
of the projection $\tXS\to X\times\MS$. Let $\mu=\mu_\tXS\in H^0(\tXS,\Omega_{\tXS/\MS})$ be the natural 1-form
on $\tXS$; it is the pull-back of the natural 1-form on $T^*X$ under the projection $\tXS\to T^*X$.

Denote by $\lS[[\lambda]]$ the pull-back of $\lS$ to $\tXS[[\lambda]]$. 
The expression $\mu+\lambda\partialS$ gives a $\lambda$-connection on $\lS[[\lambda]]$:  
\begin{equation*}
\mu+\lambda\partialS:\lS[[\lambda]]\to\lS[[\lambda]]\otimes\Omega_{\tXS[[\lambda]]/\MS[[\lambda]]}(D_r[[\lambda]]).
\end{equation*} 
So we see that over $\MS[[\lambda]]$, we have a natural family of spectral curves ($\tXS[[\lambda]]$) and line bundles
with $\lambda$-connections ($\lS[[\lambda]]$ and $\mu+\lambda\partialS$) on these curves. 
According to the generalized Theorem \ref{th:lambda-connections2}, such a family corresponds to a 
$\GL_2$-bundle $L$ on $(\MS\times X)[[\lambda]]$
equipped with a $\lambda$-connection $L\to L\otimes\Omega_X$. In other words, we obtain a morphism
$\Phi:\MS[[\lambda]]\to\Cn$ to the moduli stack of $\GL_2$-bundles with $\lambda$-connections on $X$. Clearly,
$\Phi(\MS[[\lambda]])$ is contained in $\Cnfp$ (the formal completion of $\Cn$ along $\Hgp$).

It is easy to see that Lemma \ref{lm:formal foliations}(2) applies to $\Phi:\MS[[\lambda]]\to\Cn$. Indeed,
$\MS$ is smooth, and the map 
\begin{equation*}
\lambda:\Cn\to\C:(L,\nabla,\lambda)\mapsto\lambda
\end{equation*} 
is smooth on $\Hgp\subset\Cn$; therefore,
both $\MS[[\lambda]$ and $\Cnfp$ are smooth $\lambda$-adic formal stacks. Besides, the restriction of 
$\Phi$ to $\MS\subset\MS[[\lambda]]$ (the zero locus of $\lambda$) is the natural projection
\begin{equation*}
\MS\to\MS/\zeta_0=\Hgp.
\end{equation*}

\subsection{}
Now let us verify that the linear combination $\zeta_\lambda=\zeta_0-\lambda\zeta_\infty$ equals 
$\zeta\eqd\ker(d\Phi)\subset T(\MS[[\lambda]])$ (by Remark \ref{re:combination exists}, $\zeta_\lambda$ exists
as a distribution).
As $\zeta_\lambda$ and $\zeta$ have the same rank, it suffices to check 
$\zeta_\lambda\subset\zeta$. Equivalently, for any open set $U\subset\MS[[\lambda]]$ and a vector 
field $\theta$ on $U$ that belongs to $\zeta_\lambda$, we need to check that $\theta$ belongs to $\zeta$. 

Set $U[\epsilon]\eqd U\times\spec\C[\epsilon]/(\epsilon^2)$. The vector field $\theta$ induces an automorphism 
$\Theta:U[\epsilon]\iso U[\epsilon]$ characterized by the following property:
\begin{equation*}
\Theta^*(f+\epsilon g)=f+\epsilon(g+\theta(f))\quad(\text{here }f,g\in O_U, \text{ so }f+\epsilon g\in O_{U[\epsilon]}).
\end{equation*}
We need to verify that the two compositions 
\begin{equation*}
\Phi\circ\pi,\Phi\circ\pi\circ\Theta:U[\epsilon]\to\Cn
\end{equation*}
coincide. Here $\pi:U[\epsilon]\to U\hookrightarrow\MS[[\lambda]]$ is the natural projection.

Denote by $(\tX_1,l_1,\delta_1)$ and $(\tX_2,l_2,\delta_2)$ the pull-backs of the universal family 
$(\tXS,\lS,\mu+\lambda\partialS)$ under $\pi:U[\epsilon]\to\MS[[\lambda]]$ and $\pi\circ\Theta:U[\epsilon]\to\MS[[\lambda]]$, 
respectively. Thus, $\tX_i\subset TX\times U[\epsilon]$ is a $U[\epsilon]$-family of smooth spectral curves, 
$l_i$ is a line bundle on $\tX_i$, and $\delta_i$ is a $\lambda$-connection on $l_i$ with the usual condition on the 
residues ($i=1,2$). We need to verify that $(\tX_1,l_1,\delta_1)$
and $(\tX_2,l_2,\delta_2)$ define the same $U[\epsilon]$-family of $\GL_2$-bundles with $\lambda$-connections on $X$.

According to Theorem \ref{th:lambda-connections2} (or rather its generalized version), this is equivalent to checking 
the following two statements:
\begin{enumerate}
\item The reductions of $(\tX_1,l_1,\delta_1)$ and $(\tX_2,l_2,\delta_2)$ modulo $\lambda$ coincide.
\item Let $\tX_{iu}\subset\tX_i$ be the open set where $p_{\tX_i}:\tX_i\to U[\epsilon]\times X$ is 
unramified ($i=1,2$). Then the push-forwards $(p_{\tX_1})_*((l_1,\delta_1)|_{\tX_{1u}})$ and 
$(p_{\tX_2})_*((l_2,\delta_2)|_{\tX_{2u}})$ are canonically isomorphic. Notice that the previous statement
implies $p_{\tX_1}(\tX_{1u})=p_{\tX_2}(\tX_{2u})$;
therefore, these push-forwards are $\GL_2$-bundles 
with $\lambda$-connections on the same open subset of $X\times U[\epsilon]$.
\end{enumerate}
Both statements easily follow from the definition of $\zeta_\lambda=\zeta_0-\lambda\zeta_\infty$.

This completes the proof of statements (2) and (3) of Theorem \ref{th:lambda-connections via foliations}. To prove
Theorem \ref{th:lambda-connections via foliations}(1), we need to check that the fibration $\zeta_\lambda\subset T\MS$ 
is a foliation. This translates into the vanishing of the curvature 
\begin{equation*}
\kappa:\zeta_\lambda\otimes\zeta_\lambda\to T\MS/\zeta_\lambda:\theta_1\otimes\theta_2\mapsto[\theta_1,\theta_2].
\end{equation*}
However, $\kappa$ depends on $\lambda$ algebraically, and Theorem \ref{th:lambda-connections via foliations}(2)
implies that $\kappa$ vanishes when $\lambda\in\C[[\lambda]]$ is a formal parameter. The statement follows.

\begin{remark} Theorem \ref{th:lambda-connections via foliations}(1) can be easily proved independently of
Theorem \ref{th:lambda-connections via foliations}(2). By definition, $\zeta_\lambda$ is a foliation if $\lambda=0$, 
so we can assume $\lambda\ne0$. Then $\zeta_\lambda$ can be described as an isomonodromic deformation 
(similarly to $\zeta_\infty$); the only difference is that the isomonodromic deformation uses 
twisted differential operators (where the twist depends on $\lambda$). 
\end{remark}

\section*{Acknowledgements}

I am deeply grateful to V.~Drinfeld for his constant attention to this work
and for numerous stimulating discussions. Part of this work is contained in the thesis I presented at Harvard 
University, 
and I would like to thank my readers, D.~Kazhdan and A.~Braverman, for their invaluable comments. 
D.~Kazhdan also supervised my studies at Harvard, and I would like to thank him for his guidance. 

This work owes much to the discussions I had with many mathematicians. 
I am particularly grateful to A.~Beilinson,  D.~Ben-Zvi, 
R.~Bezrukavnikov, R.~Donagi, D.~Gaitsgory, and T.~Pantev.

\bibliographystyle{plain}
\bibliography{lambdaconn}
\end{document}